\documentstyle[12pt,amssymb]{amsart}
\newmathalphabet*{\frakb}{euf}{b}{n}
\newmathalphabet*{\eusb}{eus}{b}{n}
\size{10}{18pt}\selectfont
\begin{document}
\title[One Form]
{One Form of Successive Approximation Method 
and Choice Problem}
\author{V.K.Bulitko}
\vspace*{-2.0cm}

\maketitle

\section{Introduction} 

A mathematical model of Subject behaviour choice is proposed. The 
background of the model is the concept of two preference relations
determining Subject behaviour. These are an "internal" or subjective preference 
relation and an "external" or objective
preference relation. 

The first (internal) preference relation is defined by some partial
order on a set of states of the Subject. 
The second (external) preference relation on the state set is defined by 
a mapping from the state set to another partially ordered set.
The mapping will be called evaluation mapping (function). 

We research the process of external preference
maximization in a fashion that uses the external 
preference as little as possible. On the contrary, Subject may use the
internal preference without any restriction. 

The complexity of a maximization procedure depends on the disagreement 
between these preferences.
To solve the problem we apply some kind of
the successive approximations methods. In terms of evaluation
mappings this method operates on a decomposition of
the mapping into a superposition of several standard operations 
and "easy" mappings (see the details below). Obtained in such way 
superpositions are called approximating forms. 

We construct several such forms and present two applications. One of them
is concerned with a hypothetic origin of logic. 
The other application provides a new interpretation of the
well known model of human choice by Lefebvre \cite{ForMan,HuCho}. 
The interpretation 
seems to suggest a justification different 
from the one proposed by Lefebvre himself.

\section{Scheme of Behaviour Choice Based on Two Preferences}

We consider a Subject faced with a choice among a set of states 
in the environment. Some of them may be 
better than others and some states are incomparable. 
Subject's goal is to reach
a satisfactory state (generally, a set of states). 

One fundamental feature of many real-world behaviour problems is the difference
between the evaluations of a state before and after 
the state is arrived at. We try to describe this by introducing two 
preference relations on the state set. One relation describes "internal"
system of values based on the Subject's internal representation (or model) of the world. 
The other, 
"external" relation,
is based on consequences of chosen states and reflects the actual 
nature of the interaction between the environment and the Subject. 

Unlike the objective external relation, the internal preference 
relation is intrinsic to the Subject's 
mentality for the Subject perceives the world in terms of it.
Contradictions between the internal and external preference relations 
create problems for the Subject. Besides, in general, there is a cost associated 
with obtaining information on the true external preference relation. 
Here we refer not only to the cost of accessing the information 
but also the cost of changing of Subject's behaviour
patterns. (However, we abstract from the issue of what the costs may actually be).

\medskip

Thus, informally the problem is to find a maximum in the external
preference under given restrictions on access to the information 
about the external preference. 

\medskip

However, there are no restrictions on the internal preference usage. 
Therefore, the Subject has to seek a maximum using the internal preference
as much as possible. Naturally, to get any value of using the internal preference, 
the Subject somehow needs to approximate 
the external ("leading") preference with its internal one.

It seems to us that such interpretation of the choice problem corresponds to
a certain conservatism on the side of the Subject when it is necessary 
to follow some external pressure.
Indeed, even if the Subject is aware of its incomplete and/or 
incorrect representation it often might not be able to correct it 
instantaneously. 
Therefore, it will need to refine/reconstruct its representation 
starting with what is available.

Thus, we turn to the successive approximation principle that will 
guide our further investigations.
The underlying idea is as follows. The Subject needs to follow
some part of its internal preference for as long as possible. 
Then, on the basis of accessible information on the external
preference the Subject finds the next part of the internal preference 
and uses it to proceed further. The process then repeats.
So the Subject needs a {\it scheme} to select current
parts of its intrinsic (internal) preference. 
The following section is devoted to a theory of such schemes. 

\section{An explication of successive approximation method}

Let $S$ be a set of Subject states, $(M,\le_m),(L,\le_l)$
be partially ordered sets of 
internal and external estimates correspondingly. Let $\varphi:S\to M,
\psi:S\to L$ be mappings that define internal and external preference
relations correspondingly. 

It is possible to simplify this description by introducing 
an order $\le_S$ on $S$ through the mapping $\varphi$ and 
the poset $(M,\le_m)$ as follows. Let us set 
$s\le_Ss'\iff\varphi(s)<_m\varphi(s')\vee s=s'$. 
Thus, we move from an initial description 
to the description $\langle(M,\le_m),(L,\le_l),
\psi:M\to L\rangle$ that we will use from now on.
Then $(M,\le_M)$ plays the role of $(S,\le_S)$ above
and $\psi=$ is called the evaluation mapping (function). 

It is worth noting that generally 
in each of the choice problems the Subject
gets corresponding internal and external preferences and 
evaluation mappings. These three objects can depend functionally on
some parameters of the choice problem. 
First, we consider the case of a single problem of choice. 
Later on, in the section
related to Lefebvre's model, we will generalize to a family of 
choice problems. 

If the evaluation function 
$\psi$ appears to be a monotonic mapping (i.e., the condition 
$(\forall m_1,m_2\in M)[m_1\le_mm_2\implies\psi(m_1)\le_l\psi(m_2)$]
is met), then both preference relations $(M,\le_m)$ and $(L,
\le_l)$ are compatible (coordinated) and de facto the Subject may 
follow its internal preference to reach the target state (that is, a state 
with the maximum value). 

Otherwise, it is natural to represent $\psi$ as a superposition of 
monotonic evaluation mappings from $(M,\le_m)$ to $(L,\le_l)$ and several 
connecting operations (connectives). 
We seek to obtain representations that can be used as 
instructions for successive approximations. 
In finding a representation 
of this kind such that uses as few monotonic
evaluation mappings as possible (apart from 
the connectives), we will use the external 
preference relation as little as possible. 

So our next goal is to propose some collections of 
connectives such that it will be possible 
to prove the existence of a corresponding representation with required 
features (we call it approximating form henceforth). Then, we will
demonstrate their utility for certain applications.

\subsection{Approximating forms}

The suggested version of the principle deals with 
some special but yet fairly general
representation of the evaluation function $\psi$ (operator $\psi$)
in the so-called "approximating form". 
It uses three axiomatically defined operations 
$\boxminus,\boxplus,\circledcirc$ based only on general properties
of the posets $(M,\le_m),(L,\le_l)$ as follows.

For every poset $(R,\le_r)$ the standard mappings  
$(\cdot)^{\vartriangle},(\cdot)^{\triangledown}:R\to2^R$ 
are defined by $t^{\vartriangle}=\{t'\in R|t'\le_r t\},
t^{\triangledown}=\{t'\in R|t\le_r t'\}.$

Let us suppose a binary operation  
$\boxminus:L\times L\to L$ and unary operations
$\boxplus:2^L\to L,\circledcirc:L\to L$ are defined in such a way 
that the following system $\cal A$ of axioms holds. 
\begin{description}
\item[$\cal A_1$]  
       $(\forall S\subseteq M)(\exists\breve S\subseteq S)
       [(\forall s\in S)(\exists\breve s\in\breve S)
       [\breve s\le_ms]\ \&\ (\forall\breve s,\breve s'\in\breve S)
       [\breve s\not<_m\breve s']]]$.
\item[$\cal A_2$]  
       $(\forall L',L''\subseteq L)(\forall x\in L')[(x\le_l\boxplus(L'))\&
       (L'\subseteq L''\implies\boxplus(L')\le_l\boxplus(L''))]$.
\item[$\cal A_3$]  
       $(\forall l,l'\in L)[\boxminus(l,\circledcirc(l))=l\ \&\ 
       (l\le_ll'\implies\circledcirc(l)\le_l\circledcirc(l'))]$.
\item[$\cal A_4$]
       $(\forall l,l'\in L)[l\le_ll'\implies(\exists l''\in L)
       [\boxminus(l',l'')=l\ \&\ \circledcirc(l')\le_ll'']]$.
\end{description}

For every operator $\nu:M\to L$ we call set $\frak n(\nu)=
\{(m,m')|(m\le_mm')\ \&\ (\nu(m)\not\le_l\nu(m'))\}$ {\sl non-monotonicity
domain of $\nu$}. 
If $\frak n(\nu)=\emptyset$ then $\nu$ is called monotonic operator.

{\sl{\bf Theorem 1.} Let all axioms of the system $\cal A$ be satisfied
for $(M,\le_m),(L,\le_l)$ and lengths of all increasing chains in
$(M,\le_m)$ not exceed some integer $D$.
Then for every $\psi:M\to L$ there exists a representation
$\psi=\boxminus(\varphi_1,
\boxminus(\varphi_2,\boxminus(\varphi_3,\dots)))$ that all 
$\varphi_i,i=1,2,3\dots$ are monotonic mappings from  
$(M,\le_m)$ to $(L,\le_l)$.

Furthermore, the number of occurrences of the operation $\boxminus$ in this 
representation does not exceed $D$. 
}

{\sl\bf Proof.} 
Let us reduce the problem for a given operator    
$\psi$ to the same problem for a simpler operator
$\psi_1$ such that the following holds
$\psi=\boxminus(\varphi_1,\psi_1)$ and $\frak n(\psi_1)
\subsetneqq\frak n(\psi)$.

First, we define $M_1=\{x\in M|\frak n(\psi)\cap(x^{\vartriangle}
\times x^{\vartriangle})\neq\emptyset\}$,
$M^1=\overline{M_1}$ ¨ 
\begin{center}
$\varphi_1(x) =
\begin{cases}
\boxplus(\psi(x^{\vartriangle})),&  x\in M_1,\\
\psi(x),&   x\in M^1.
\end{cases}$
\end{center}
Then we set $\psi_1(x)$ to any such $z\in L$ that
$\boxminus(\varphi_1(x),z)=\psi(x)\ \&\ 
\circledcirc(\varphi_1(x))\le_lz$ if
$\varphi_1(x)\neq\psi(x)$. Otherwise, we set 
$\psi_1(x)=\circledcirc(\psi(x))$.

Existence of the element $z$
in the definition is guaranteed by axioms $\cal A_3,\cal A_4$.
Now, the equality $\psi(x)=\boxminus(\varphi_1(x),\psi_1(x))$ holds
because of the definitions of $\varphi_1,\psi_1$. 

Let us prove that operator $\varphi_1:(M\le_m)\to(L,\le_l)$ is
a monotonic one. 

First,    
$\varphi_1=\psi$ on $M^1$ and we may use the condition
$x,y\in M^1\& x\le_my\implies\psi(x)\le_l\psi(y)$. 
Indeed, otherwise 
$\psi(x)\not\le_l\psi(y),x\le_my,\psi(x)\neq\psi(y)$ and, therefore,
$y\in M_1\cap M^1$. However, $M^1\cap M_1=\emptyset$ 
which leads to a contradiction.  

Second, $\varphi_1$ maps $(M_1,\le_m)$ into
$(L,\le_l)$ monotonically in accordance with $\cal A_2$. 

Finally, let us consider the "mixed" case when $x\in M^1,y\in M_1$
and all elements of $M$ are comparable with respect to $\le_m$. 
It is clear that $y\le_mx$ is impossible since condition
$z\in M_1\implies z^{\triangledown}\subseteq M_1$ immediately follows 
from the definition of $M_1$.

Thus, it remains to consider the possibility of $x\le_my$. 
In such a case
$\varphi_1(y)=\boxplus(\psi(y^{\vartriangle}))\ge_l\psi(x)$
in accordance to $\cal A_2$. On the other hand, $\psi(x)=\varphi_1(x)$
on $M^1$ follows from the definition of $\varphi_1$. 
Hence, operator $\varphi_1$ is monotonic. 

We are now ready to prove the last assertion of the theorem. For that it is sufficient
to demonstrate the inclusion $M^1\cup\breve{M_1}\subseteq M^2$.
Here  $M^2,M_2$ 
are defined for $\psi_1$ in the same way as $M^1,M_1$ were defined
for $\psi$ above. $\breve{M_1}$ is the set of all minimal elements 
of set $M_1$, see $\cal A_1$. 
Namely: $M^2=\overline{M_2}$ and 
$M_2=\{x\in M|\frak n(\psi_1)\cap(x^{\vartriangle}
\times x^{\vartriangle})\neq\emptyset\}$. 

We now have 
$M_2\subseteq(M_1\setminus\breve{M_1})$ and   
$\frak n(\psi_1)\subseteq\frak n(\psi)\setminus\breve M_1\times M_1$. 
So the sequence $M_1\supsetneqq M_2\supsetneqq M_3\supsetneqq\dots$ 
ends at a step with the number that can not 
be greater than the highest
of lengths of the increasing chains in poset $(M,\le_m)$. Indeed, since
$\breve M_2\subseteq M_1\setminus\breve M_1$ then in accordance
with $\cal A_1$
for every element $y\in\breve M_2$ there exists some $x\in\breve M_1$
such that $x<_ly$. Therefore, one can choose an increasing chain 
of representatives of sets $\breve M_1,\breve M_2,\breve M_3,\dots$
which are mutually disjoint sets. 

We will now prove that $M^1\cup\breve{M_1}\subseteq M^2$.
First, $\varphi_1(x)=\psi(x)$ holds for every $x\in M^1$.
From here $\psi_1(x)=\circledcirc(\psi(x))$.
However, mapping $\psi_1$ is 
monotonic on $M^1$ in view of $\cal A_3$ and
since $\psi$ is monotonic on $M^1$. 
So $(M^1\times M^1)\cap\frak n(\psi_1)=\emptyset$ and therefore 
$M^1\subseteq M^2$.

Further, let $x,y\in M^1\cup\breve{M_1}$ and $x\le_m y$. Then
we can show that $\psi_1(x)\le_l\psi_1(y)$. 
Indeed, the case $x,y\in M^1$ was considered above. 
The case $x,y\in\breve{M_1}$ is impossible since all elements of 
$\breve{M_1}$ are incomparable by the definition.
Above, we saw that $x\in M_1\ \&\ x\le_my\implies y\in M_1$. 
Besides $M^1\cap M_1=\emptyset$. 
Therefore, $x\in M^1,y\in\breve M_1$ is the only case remaining 
to consider. 
By definition $\psi_1(x)=\circledcirc(\psi(x))$ and 
relation $\boxminus(\varphi_1(y),\psi_1(y))=\psi(y)$ holds. Moreover, 
$\psi(y)<_l\varphi(y)$. In accordance with $\cal A_4$ we have 
$\circledcirc(\varphi_1(y))\le_y\psi_1(y)$. Hence,
$\psi_1(x)\le_l\psi_1(y)$ takes place since $\circledcirc$ is a monotonic
operation in view of $\cal A_3$ and 
$\psi(z)\le_l\varphi_1(z),z\in M$ in accordance to $\cal A_2$ and the 
construction. \ \ \ $\Box$

Let us denote by $\cal M$ the class of all monotonic mappings from
$(M,\le_M)$ to $(L,\linebreak\le_L)$. Also let $S(\varphi)=
\{x|\varphi(x)>_L\circledcirc(x)\},\varphi\in\cal M$. 

{\sl{\bf Corollary 1.} Under the conditions of theorem 1 for every  
$\psi:M\to L$ there exists a substitution
$p:\{z_1,\dots,z_{D+1}\}\to\cal M$ such that 
$S(p(z_n))\subseteq S(p(z_{n+1})),
n=\overline{1,D},$ and 
$\psi=\bold{Sb}^{z_1\ \ \dots\ \ z_{D+1}}_{p(z_1)\dots p(z_{D+1})}
\boxminus(z_1,\boxminus(z_2,\boxminus(\dots\boxminus(z_D,z_{D+1})\dots)))$. 
}

{\sl\bf Proof.} 
Let us fix a formula $\Phi(z_1,\dots,z_{D+1})=
\boxminus(z_1,\boxminus(z_2,(\dots
\boxminus(z_D,z_{D+1})\dots)))$ and 
consider substitutions of monotonic functions instead of variables
$z_1,\dots,z_{D+1}$  when their results are determined.

According to theorem 1 for every
$\psi:M\to L$ there exists representation
$\psi=\boxminus(\varphi_1,
\boxminus(\varphi_2,\boxminus(\varphi_3,\dots)))$ where all 
$\varphi_i,i=1,2,3\dots$ are monotonic mappings from 
$(M,\le_M)$ into $(L,\le_L)$. The number of occurrences of operation
$\boxminus$ in this representation does not exceed $D$. 

Now, if for a given mapping $\psi$ we obtain representation
$\psi=\boxminus(\varphi_1,\boxminus(\varphi_2,\boxminus(\varphi_3,\linebreak
\dots\boxminus(\varphi_{k},\varphi_{k+1}))))$ 
with $k<D$, then one can always continue
the expression on the right side of the representation till 
$\boxminus(\varphi_1,\boxminus(\varphi_2,\boxminus(\varphi_3,
\dots)\varphi_{D+1}))$. For that it is sufficient to set
$$
\varphi_{i}(x)=\circledcirc(\varphi_{i-1}(x)),i=\overline{k+2,D+1}.
$$
In accordance with axiom $\cal A_3$ the obtained functions
are monotonic and 
$$
\psi=\boxminus(\varphi_1,\boxminus(\varphi_2,\boxminus(\varphi_3,
\dots\boxminus(\varphi_{D},\varphi_{D+1}))))
$$ 
$\ \ \Box$

For the application below we will need some special monotonic
functions in approximating forms. To define the functions we come up with 
the  following auxiliary construction.
Based on axiom $\cal A_1$ let us split set $M$:\\ 
$M_1=M^{\bot}$;\\
$M_{n+1}=(M\setminus\underset{j\le n}{\cup}M_j)^{\bot}$.\\
Here $M_j$ consists of the elements that are incomparable in $(M,\le_M)$ 
for any $j$.

We denote by $\theta$-function of rank $i$ any monotonic mapping
$\theta:M\to L$ such that $\theta(x)=\circledcirc(x)$ for all 
$x\in\underset{j<i}{\cup}M_j$ as well as $\theta(x)\in\max(L,\le_L)$ 
for all $x\in\underset{j>i}{\cup}M_j$. The rank of a given function
$\theta$ is denoted as $\rho(\theta)$. Let $\Theta$ be the class of all
$\theta$-functions.

{\sl{\bf Corollary 2.} Let conditions of theorem 1 be fulfilled, 
$D$ be the exact upper bound of lengths of increasing chains in
$(M,\le_M)$ and $(L,\le_L)$ contain the greatest element $\gamma$. 
Then for any mapping $\psi:M\to L$ 
there exists a substitution $p:\{z_1,\dots,z_{D+1}\}\to\Theta^{\mu}$, 
such that  $\rho(p(z_i))=i,i=\overline{1,D+1}$ and 
$\psi=\bold{Sb}^{z_1\ \ \dots\ \ z_{D+1}}_{p(z_1)\dots p(z_{D+1})}
\boxminus(z_1,\boxminus(z_2,\boxminus(\dots\boxminus(z_D,z_{D+1})\dots)))$.
}

{\sl\bf Proof.} 
We use induction on $D$. In the case of $D=0$ the statement is obvious
since there are no restrictions for $\theta$-functions.
Therefore, $\psi\in\Theta$. 

{\sl Induction step:}
Let us define $\theta_1$ of rank 1 in the following manner:
\begin{center}
$\theta_1(x) =
\begin{cases}
\psi(x),&  \text{ if }x\in M_1,\\
\gamma,&   \text{ else}.
\end{cases}$
\end{center}
Then we may write $\psi=\boxminus(\theta_1,\psi_1)$ where for
$\psi_1$ we have $\psi_1(x)=\circledcirc(x)$ if $x\in M_1$ else
$\psi_1(x)$ satisfies $\boxminus(\gamma(x),\psi_1(x))
=\psi(x)$. 

It remains to obtain the desirable representation of $\psi_1$ on the set
$M\setminus M_1$ with the partial order induced by $\le_M$.
Since the length of the longest increasing chain in
$(M\setminus M_1,\le_M)$ is $D-1$ we may use the induction supposition.
$\ \ \ \Box$

{\bf Remark.} 
Length of the representation obtained in theorem 1 can
be essentially lower than lengths of the representations suggested by
the corollaries.
From our initial standpoint the lower the length is the better.
However, sometimes we will need special forms of approximating functions.

Let us suppose a binary operation  
$\boxminus,\uplus:L\times L\to L$ and unary operations
$\circledcirc:L\to L$ are defined in such a way 
that the system $\cal B_1\ -\ \cal B_4$ of axioms takes place. 
Here $\cal B_i$ coincides with $\cal A_i$ for $i=1,3,4$. Also
\begin{description}
\item[$\cal B_2$]  
       $(\forall x,y\in L)[x,y\le_l\uplus(x,y)]$.
\end{description}

{\sl{\bf Theorem 2.} Let all axioms of the system $\cal B$ be satisfied
for $(M,\le_m),(L,\le_l)$ and $(M,\le_m)$ 
and lengths of all increasing chains in
$(M,\le_m)$ do not exceed some integer $D$.
Then for every $\psi:M\to L$ there exists a representation
$\psi=\boxminus(\varphi_1,
\boxminus(\varphi_2,\boxminus(\varphi_3,\dots)))$ where all 
$\varphi_i,i=1,2,3\dots$ are monotonic mappings from  
$(M,\le_m)$ to $(L,\le_l)$.

The number of occurrences of the operation $\boxminus$ in this 
representation does not exceed $D$. 
}

{\sl\bf Proof.} Firstly, in the case when 
$(\forall x\in M)[|x^{\vartriangle}|<\infty$ is true
we can prove this theorem using theorem 1. 
For that we will only need to note that in this case
it is possible to replace $\boxplus(\psi(x^{\vartriangle})$ with 
any expression of the kind $\uplus(psi(z_1),\uplus(\dots
\uplus(\psi(z_{n_1},\psi(z_{n})\dots))$. Here $z_1,\dots,z_n$ is an
enumeration of the finite set $x^{\vartriangle}$. 
Indeed, in the proof of theorem 1 we used axiom $\cal A_2$ 
only for subsets of $L$ of the form $\psi(x^{\vartriangle})$. 
Thus, it is sufficient to check that axiom $\cal A_2$ is
met for sets of the kind $\psi(x^{\vartriangle})$. 
This check is trivial on the basis of axiom $\cal B_2$ 
for operation $\uplus$. 

Otherwise, when there are infinite sets $x^{\vartriangle}$ we can 
make use of the same scheme for the connective
$\boxplus$ based on the condition of increasing chain finiteness in 
$(M,\le_m)$. For that let us enumerate elements 
$z_1,z_2,\dots,z_n,\dots$ of set 
$x^{\vartriangle}$ for a given $x\in M$. In parallel we will enumerate
expressions 
$\uplus(\psi(z_1),\psi(z_2)),\ 
\uplus(\uplus(\psi(z_1),\psi(z_2)),\linebreak
\psi(z_2)),
\dots,\uplus(\uplus(\dots\uplus(\uplus(\psi(z_1),\psi(z_2)),
\psi(z_2))\dots),z_n),\dots$

By axiom $\cal B_2$ the values of these expressions do not decrease in
$(L,\le_L)$. In view of the finiteness supposition for increasing chains in
$(L,\le_L)$ the sequence of computed values becomes stable from
some place. We set $\varphi_1(x)$ to this final value. 

Thus, $\varphi_1$ is monotonic mapping and
$\psi(x)\le_L\varphi_1(x),x\in M$. The last part of the proof
is analogous to the corresponding part of theorem 1. $\ \ \Box$

In place of or together with $\cal A$ the dual axiom system
$\cal A^{\star}$ can be fulfilled. This can be shown by replacing
$\le$ with $\ge$ and $\boxplus,\boxminus,\circledcirc$ with
$\boxplus^{\star},\boxminus^{\star},\circledcirc^{\star}$ correspondingly:
\begin{description}
\item[$\cal A^{\star}_1$]  
       $(\forall S\subseteq M)(\exists\breve S\subseteq S)
       [(\forall s\in S)(\exists\breve s\in\breve S)
       [\breve s\ge_ms]\ \&\ (\forall\breve s,\breve s'\in\breve S)
       [\breve s\not<_m\breve s']]]$.
\item[$\cal A^{\star}_2$]  
       $(\forall L',L''\subseteq L)(\forall x\in L')[(x\ge_l
       \boxplus^{\star}(L'))\&
       (L'\subseteq L''\implies\boxplus^{\star}(L')\ge_l
       \boxplus^{\star}(L''))]$.
\item[$\cal A^{\star}_3$]  
       $(\forall l,l'\in L)[\boxminus^{\star}(\circledcirc^{\star}(l),l)=l
       \ \&\ 
       (l\ge_ll'\implies\circledcirc^{\star}(l)\ge_l
       \circledcirc^{\star}(l'))]$.
\item[$\cal A^{\star}_4$]
       $(\forall l,l'\in L)[l\ge_ll'\implies(\exists l''\in L)
       [\boxminus^{\star}(l'',l')=l\ \&\ 
       \circledcirc^{\star}(l')\ge_ll'']]$.
\end{description}
Then the dual theorem holds:

{\sl{\bf Theorem 1$^{\star}$.} 
Let all axioms of the system $\cal A^{\star}$ 
be fulfilled for posets $(M,\le_m),(L,\le_l)$, 
operators $\boxplus^{\star},\boxminus^{\star},\circledcirc^{\star}$;
and the lengths of all increasing chains in
$(M,\le_m)$ do not exceed some integer $D$.
Then for every operator $\psi:M\to L$ there exists representation
$\psi=\boxminus^{\star}(\dots\boxminus^{\star}(\boxminus^{\star}
(\varphi_{n+1},\varphi_{n}),
\varphi_{n-1})\dots,\varphi_1)$ where all
$\varphi_i,i=1,\dots,n,n+1,$ are monotonic mappings from
$(M,\le_m)$ to $(L,\le_l)$.

The number of occurrences of the operation $\boxminus$ in this 
representation does not exceed $D$. 
}

The dual theorem is related to the dual axiom system 
$\cal B^{\star}$.

\begin{description}
\item[$\cal B^{\star}_1$]  
       $(\forall S\subseteq M)(\exists\breve S\subseteq S)
       [(\forall s\in S)(\exists\breve s\in\breve S)
       [\breve s\ge_ms]\ \&\ (\forall\breve s,\breve s'\in\breve S)
       [\breve s\not<_m\breve s']]]$.
\item[$\cal B^{\star}_2$]  
       $(\forall x,y\in L)[x,y\ge_l\uplus^{\star}(x,y)]$.
\item[$\cal B^{\star}_3$]  
       $(\forall l,l'\in L)[\boxminus^{\star}(\circledcirc^{\star}(l),l)=l
       \ \&\ 
       (l\ge_ll'\implies\circledcirc^{\star}(l)\ge_l
       \circledcirc^{\star}(l'))]$.
\item[$\cal B^{\star}_4$]
       $(\forall l,l'\in L)[l\ge_ll'\implies(\exists l''\in L)
       [\boxminus^{\star}(l'',l')=l\ \&\ 
       \circledcirc^{\star}(l')\ge_ll'']]$.
\end{description}
Then the dual theorem holds:

{\sl{\bf Theorem 2$^{\star}$.} Let all axioms of the system $\cal A^{\star}$ 
be fulfilled for posets $(M,\le_m),(L,\le_l)$, 
operators $\uplus^{\star},\boxminus^{\star},\circledcirc^{\star}$
and  
the lengths of all increasing chains in
$(M,\le_m)$ do not exceed some integer $D$.
Then for every operator $\psi:M\to L$ there exists representation
$\psi=\boxminus^{\star}(\dots\boxminus^{\star}(\boxminus^{\star}
(\varphi_{n+1},\varphi_{n}),
\varphi_{n-1})\dots,\varphi_1)$ where all
$\varphi_i,i=1,\dots,n,n+1,$ are monotonic mappings from
$(M,\le_m)$ to $(L,\le_l)$.

The number of occurrences of the operation $\boxminus$ in this 
representation does not exceed $D$. 
}

Below we refer to all these representations as approximating forms.

\section{Possible Origin of Logic}

It is easy to arrive at the classical two-valued propositional logic now. 
For that it is sufficient to choose
$(\{0,1\},0\le1)$ as $(L,\le_l)$ and 
the standard poset $(\cal B^n,
\preccurlyeq)$ on boolean cube $\cal B^n$ as poset $(M,\le_m)$.
It is well known that every finite poset can be isotonically 
included into $(\cal B^n,\preccurlyeq)$ for the appropriate $n$. 

It is also well known that poset $(\cal B^n,\preccurlyeq)$
is a self-dual poset for any  $n$. Therefore, both 
above introduced representations
take place in this case. 

{\sl{\bf Lemma.} 1) The system of posets $(\cal B^n,\preccurlyeq),
(\cal B,\le)$ as $(M,\le_m),(L\le_l)$ correspondingly and operation
$\to$ as $\boxminus^{\star}$, operation $\Bbb{\bold1}:\cal B^n\to\{1\}$ as
$\circledcirc^{\star}$, and operation  
$\underset{\vec{\beta}\preccurlyeq\vec{\alpha}}{\&}\vec{\alpha}$
as $\boxplus^{\star}(\beta^{\triangledown})$ fulfill the axiom
set $\cal A^{\star}$.\\
2) The system of posets $(\cal B^n,\preccurlyeq),
(\cal B,\le)$ as $(M,\le_m),(L\le_l)$ correspondingly and operation
$\to$ as $\boxminus^{\star}$, operation $\Bbb{\bold1}:\cal B^n\to\{1\}$ as
$\circledcirc^{\star}$, and operation $\&$
as $\uplus^{\star}$ obeys the axiom set ${\cal B^+}^{\star}$. 
}

{\sl\bf Proof.} This can be shown via a routine check of the axioms.

The direct corollary of this lemma and  theorems above is the following

{\sl{\bf Theorem 3.} In the special case  
of finite "internal" orders $(M,\le_m)$
and linear "external" orders $(L,\le_l),|L|=2$,  
approximating forms from each of theorems 1 and 2 and their dual ones 
generate all
formulae of the classical propositional logic (within logical equivalence). 
}

As a result, this interesting statement follows.

{\sl{\bf Corollary 3.} Every $n$-argument logical (boolean) function 
$f$ can be represented by the implicative normal form 
$f=P_k\to P_{k-1}\to\dots\to P_1$, where 
$k\le n$, and $P_i,i=\overline{1,k},$ are monotonic boolean function.
}

It is remarkable that just the dual approximating forms present the usual 
propositional implication. One may then wonder why the operation $\to^{\star}$ is not 
present 
in natural languages? In our opinion, the main reason 
is that the dual approximating forms of theorems 1$^{\star}$, 2$^{\star}$ 
begin with 
a given operator $\psi$ and approximate it by means of successive 
simplifications: $\psi_1=\boxminus^{\star}(\psi,\varphi_1),\psi_2=
\boxminus^{\star}(\psi_1,\varphi_2),\dots$
while $\psi_i$ is not a monotonic operator (i.e., not an "easy" one). 
Thus, the approximation begins with a target unlike in  
the case of the approximating forms in theorems 1,2. 

Now one can consider the classical two-valued
propositional logic merely as a realization of 
the above-mentioned principle of successive approximations for 
the problem of decision-making 
within Subject-environment survival framework. 

Thus, from this viewpoint, 
the classical propositional logic can take its beginning 
from the survival problem. It is also important that
this hypothetical origin of logic appears
quite natural. 

\section{What stands behind Lefebvre's model}

Lefebvre suggested 
a model of Subject facing a choice of an alternative out of a set. 
In his model the Subject is represented by 
the function $X_1=f(x_1,x_2,x_3)$ where 
$X_1,x_1,x_2,x_3$ run over the $[0,1]$ segment. 
As \cite{ForMan, HuCho} presents it:
the value of $X_1$ is interpreted as the readiness to choose 
a positive pole 
with probability $X_1$, and the value of $x_3$ - as the Subject's 
plan or intention to choose a positive pole with probability $x_3$.
Variables $x_1$ and $x_2$ represent the world influence on the subject.

This function $f$  is required to obey
the following axioms introduced by Lefebvre:
\begin{description}
\item[$\cal L_1$] 
$(\forall x_3\in[0,1])(f(0,0,x_3)=x_3)$ - "the axiom of free choice";
\item[$\cal L_2$] 
$(\forall x_3\in[0,1])(f(0,1,x_3)=0)$ - "the axiom of credulity";
\item[$\cal L_3$] 
$(\forall x_2,x_3\in[0,1])(f(1,x_2,x_3)=1)$ - 
"the axiom of non-evil-inclinations";
\item[$\cal L_4$] 
$(\forall i,j,k)[\{i,j,k\}=\{1,2,3\})\implies
(\forall x_j,x_k\in[0,1])(\exists c,c'\in\Bbb R)
(\forall x_i\in[0,1])[f(x_1,x_2,x_3)=cx_i+c']]$ - "the postulate of simplicity".
\end{description}

By means of the model Lefebvre gave explanations of several psychological
experiments thusly putting his model under the spotlight 
(e.g., see bibliography in \cite{HuCho}). 

The following question is still open: Is the model only a compact 
representation (i.e., a "roll-up") of certain experimental data 
or it describes some fundamental structure governing human behavior?  

In order to substantiate his model, Lefebvre used, in particular, known
"anthropic principle" \cite{HuCho}.
In our opinion, the justification presented by Lefebvre while being 
appealing does not appear 
entirely sound and bullet-proof. The specific comments are presented 
in \cite{Bul}. 
In the following we suggest an alternative justification for the model. 
Namely, we develop the 
approach mentioned in \cite{Bul} using the above constructed theory 
of approximating forms.

First, we show how it is possible to eliminate "the postulate of simplicity"
introducing the notion of a pure L-ensemble. The last concept reduces the general
case to the boolean case. This step leads to the boolean order for the
external preference relation $(L,\le_l), L=\{0,1\}$. 
Second, we will show that the system of the first
three axioms by Lefebvre 
can be replaced with a postulate of special poset $(M,\le_m)$.
Namely, this poset can be chosen in the
form of a linear ordered three-element set. We suggest a natural 
interpretation of this form of poset $(L,\le_l)$. 
Then Lefebvre's function $f$ follows from one of our approximating forms. 

\subsection{Lefebvre's ensembles}

It is easy to check that in the boolean case $X_1,x_1,x_2,$ $x_3\in\{0,1\}$
the axioms $\cal L_1-\cal L_3$ completely define $f$. Namely,
in this case $f(x_1,x_2,x_3)=(x_3\to x_2)\to x_1$. 
(The "postulate of simplicity" $\cal L_4$ sets $f$ on the interior
of the three-dimensional cube $[0,1]^3$ in the real-valued case. 
A methodological
criticism of the postulate is expounded in \cite{Bul}).

Let us consider a set $Q$ of Subjects $s_i$ with each 
being described by the probabilistic collection $\tilde{\alpha}_i$ of values of
the boolean variables $(n_1,n_2,n_3)$. 
Let us assume that the probability of encountering 
a Subject with a collection $\tilde{\alpha}$ of the variable values in $Q$
is equal to $p_{\tilde{\alpha}}$. 

If behavior $z_i$ each $s_i\in Q$ is described with the function
$n_3\to n_2\to n_1$ then we refer to $Q$ as {\it 
the Lefebvre's ensemble ($L$-ensemble or {\rm simply} ensemble) 
$\langle Q,P\rangle$ with characteristic $P=(p_0,\dots,p_7)$}. 
Besides, we call elements of the $L$-ensemble $L$-Subjects.   
(Here  $p_k$ denotes   
$p_{\tilde{\alpha}}$ and  $k$ is the decimal representation of the binary 
sequence $\tilde{\alpha}$). 

Ensemble $\langle Q,P\rangle$ averaging Boolean variables
$n_1,n_2,n_3,z_i$ yields real numbers 
$x_1,x_2,x_3,z\in[0,1]$. Given the truth table of the Boolean function
$n_3\to n_2\to n_1$ elementary probabilistic considerations lead to
the following equalities:
\begin{eqnarray}
1=\overset{7}{\underset{k=0}{\Sigma}}p_k,\\
x_1=p_4+p_5+p_6+p_7,\\
x_2=p_2+p_3+p_6+p_7,\\
x_3=p_1+p_3+p_5+p_7,\\
z=p_1+p_4+p_5+p_6+p_7.
\end{eqnarray}
It is therefore reasonable to ask for which $L$-ensembles 
$\langle Q,P\rangle$ values of $x_1,x_2,x_3,z$ satisfy 
Lefebvre's equation $z=x_1+(1-x_1-x_2+x_2x_3)x_3$. 
 
The following examples show that, generally speaking, $z\neq f(x_1,x_2,x_3)$.
Indeed, let us set $p_1=p_2=p_3=p_4=p_5=p_6=p_7=0,1$. Then 
$x_1=x_2=x_3=0.4$ and $f(x_1,x_2,x_3)=0.544$. However, the ensemble average
$z$ equals $0.5$. 
Interestingly enough, the difference can be quite substantial
as the following example demonstrates.
Namely, $p_0=p_1=p_2=p_4=p_6=p_7=0,p_3=p_5=0.5$ correspond to 
$x_1=x_2=0.5,x_3=1$. Then $z=0.5$ but $f(0.5,0.5,1)=0.75$. Thus, the error
is at least 30\%. 

On the other hand, the equality $z=f(x_1,x_2,x_3)$ is met for all
possible (i.e., obeying equations (1)-(4)) characteristics $P$ when
$(x_1,x_2,x_3)\in\{(x_1,x_2,x_3)|x_1=1\}\cup\{(x_1,x_2,
x_3)|x_2=0\}\cup\{(x_1,x_2,x_3)|x_2=1\}\cup\{(x_1,x_2,x_3)|x_3=0\}$. 

{\sl{\bf Theorem 4.} For every collection $x_1,x_2,x_3\in[0,1]$ 
there exists $L$-ensemble $\langle Q,P(x_1,x_2,x_3)\rangle$ with 
  characteristic $P(x_1,x_2,x_3)$ such that $z=f(x_1,x_2,x_3)$. 
}

{\sl\bf Proof.}
Let us consider three independent Boolean random 
variable $\zeta,\eta,\theta:\Bbb N\to\{0,1\}$
with the mean values $x_1,x_2,x_3$ correspondingly. 
Then 
random variable $(\zeta,\eta,\theta):\Bbb N\to\{0,1\}^3$ 
runs over the desired ensemble $\langle Q,P(x_1,x_2,x_3)\rangle$. 
For the $i$-th component of the characteristic $p_i(x_1,x_2,x_3)=
\underset{j=1,2,3}{\Pi}(1-\sigma_j+(-1)^{1-\sigma_j}x_j)$ is true where
$i=\underset{j=1,2,3}{\Sigma}2^{\sigma_j}$. 
The verification by substitution shows that 
the interrelations (1)-(4) are fulfilled and if $z$ satisfies (5), 
then $z=f(x_1,x_2,x_3)$. 
$\ \ \Box$

We call the ensembles described in this theorem 
{\it pure Lefebvre's ensembles ($PL$-ensembles)}. 
Thus, a $PL$-ensemble is a collection of $L$-Subjects with random parameters
$(n_1,n_2,n_3)$ distributed independently in such a way
that the probability $\bold P\{n_i=1\}$ equals the given
number $x_i\in[0,1],i=1,2,3$.

The descriptions of behaviour
constructed by means of $L$-ensembles can be thinner 
than the descriptions "smoothed"
by using Lefebvre's function $f$ for some aspects. 
For example, let us consider how "golden section"
for categorization of stimuli without measurable intensity
can be explained in terms of Lefebvre's theory (\cite{HuCho}, p.51) 
and in terms of $PL$-ensembles. 

In this case Lefebvre completes his "Realist condition" 
$x_3=f(x_1,x_2,x_3)$ with equations $x_1=x_2,x_1=1-x_3$. (A justification 
is given in \cite{HuCho}, p.51). In turn, that yields the
equation $x_3^3-2x_3+1=0$ for the choice of $x_3$. 
One possible solution is the well known "golden section" $x_3=\frac{\sqrt 5-1}{2}$. 

Following the alternative approach suggested in this paper, we construct 
the desired $PL$-ensemble  
by first postulating the Boolean
"Realist condition" $n_3\to n_2\to n_1=n_3$. Then considering
the truth area $R=\{000,001,010,101,111\}$ of the condition we form 
the ensemble with the help of Boolean random variables $\zeta,\eta,\theta$
in the following fashion. The variables $\zeta,\eta$ are independent
with the mean value $1-x_3$, and the value of the random variable $\theta$ 
depends on the values of $\zeta,\eta$ in accordance with the table:
\begin{center}
\begin{tabular}{|c|c|c|}
\hline 
$\zeta$ & $\eta$ & $\theta$ \\
\hline
0 & 0 & 0,1\\
\hline
0 & 1 & 0 \\
\hline
1 & 0 & 1  \\
\hline
1 & 1 & 1 \\
\hline
\end{tabular}.
\end{center}
It is important that in the first line of the table value 1 is chosen with
the probability of $x_3$. Then if $x_3$ satisfies $x_3^3-2x_3+1=0$ then we obtain 
the desired $PL$-ensemble. 
Indeed, every element of the ensemble is
a "Realist" and the probability to encounter an $L$-Subject with parameters
$(n_1,n_2,1)$ is determined by solutions to
the equation $x_3^3-2x_3+1=0$.
Finally, we arrive at the "golden section" choosing the corresponding  
solution exactly as it was done by Lefebvre.
  
We believe that the $L$-ensemble tool introduced in this paper 
opens new opportunities
for Lefebvre's theory and its applications. Indeed, 
the ensemble structure is a new powerful parameter for modeling. 
It is possible to explain some deviations of the actual values 
of variable $X_1$ in real-world experiments by means of
the corresponding deviations of the real $L$-ensembles 
from the $PL$-ensembles.
Thus, the dynamics of this parameter open a new research avenue. 

\subsection{Application of approximating forms}

We will now show how one can arrive at Lefebvre's model on the basis of
the theory presented earlier in this paper.

First, we determine appropriate internal and external preferences. 
Because of the binary choice in Lefebvre model it is naturally to 
take $(\{0,1\},\le)$ as the external order. (Here $\le$ is the usual
order on the set of integers). 

Second, according to the interpretation of variables $x_1,x_2,x_3$
given by Lefebvre, the values of these variables describe directions of 
impulses (motivations) pushing the Subject to the positive or the negative pole. 
Indeed, $x_1$ corresponds to an impulse exerted by the external world,
$x_2$ corresponds to an impulse exerted by Subject's experience, and, 
finally, $x_3$ corresponds to Subject's will. 

So on one hand, $x_1,x_2,x_3$ are connected to the motivations. 
On the other hand, at any decision node these variables  
have boolean values. 
Furthermore, the choice of some of these values represents
the result of the decision node. 

In our approach these two sides of variables work
simultaneously. We describe the impulses ("pure motivations")  
by partial orders (whereas 
results of Subject's choice are numbers 0 or 1). Two possible values of
a variable present two possible pure motivations for this variable.
Our choice of domain of
these partial orders is based on the following reasons. 

These six (two specific pure motivations for every variable of 
$\{x_1,x_2,x_3\}$) 
partial orders are basic and their interaction would
determine Subject's choices within our frame of 
two-preference decision-making. The decision making is done in two stages.
At the first stage some of the given pure motivations (i.e., some variables)
are chosen. At the second stage the Subject proceeds to the pole
associated by Lefebvre's 
interpretation with the given value of the variable. This means that
Lefebvre's state set $M_*$ has to be $\{x_1,x_2,x_3\}$. 
In our scheme the chosen state has to maximize external value that 
is computed with the current evaluation mapping $\psi$. 
Hence, $\psi$ sets the external preference and the latter, in turn,
determines Subject's choices (decisions).

We now seem to come to the conclusion that it is the {\it interaction} of 
pure motivations 
that produces these external preferences or equivalently $\psi$. So 
the external preferences have to be some sort of "mixture" of
pure motivations. 
(Here one can notice a vague analogy with quantum mechanics.) 

Maximizing Subject's adaptation abilities leads to the best survival
chances. Therefore, one seeks a universal "mixing" procedure. 
Corollary 2 tells us that such a procedure can be attained using 
the universality of the 
corresponding approximating form (in our case $D=2$ because
two is the upper bound of lengths of the longest increasing chains 
possible in posets of three elements): 	
$$
\psi=\bold{Sb}^{z_1,z_2,z_3}_{p(z_1),p(z_2),p(z_3)}
\boxminus(z_1,\boxminus(z_2,z_3)).
$$
Here $p$ runs over the class of special functions $\Theta$. Given 
the chosen external order we can set $\boxminus=\nrightarrow$ (see
lemma 1 above. $\nrightarrow$ is the connective 
dual for implication). 
Hence it follows that $\Theta$ is the set of pure 
motivations (impulses) in this case. So $|\Theta|=6$. 
Therefore, the sought order on $M_*$ has to be a linear.

It is natural to deem that this internal linear poset 
reflects the common division of the time axis in 
three periods: "past", "present" and "future".
Then the current representation of the world (variable $x_2$
in Levebvre's model) corresponds to the point "present" and the Subject's intention
($x_3$) corresponds to the point "future". Thus, the remaining variable $x_1$ ought to correspond
to the point "past". Such assignment
appears natural because the pressure put on the Subject by the environment
is the background of the decision-making problem itself. 
Thus, we arrive at the internal preference relation $(M_*,\le_*)$, where
$M_*=\{x_1,x_2,x_3\},x_1<_*x_2<_*x_3$.

It may seem that $x_2<_*x_1$ ought to hold since we interpret $x_2$ as the "past experience" 
and $x_1$ as the "current pressure of the environment".

However, we should keep in mind that we are currently 
dealing with the {\it internal order} on
states in the process of decision making. In that process "past
experience" $x_2$ serves a role of Subject's "current base" and it is 
$x_1$ that initiates
the decision making. $x_3$ is merely a means to produce a solution and as such 
is most likely
related to the future. (Note that such crude models often cover several 
various factors with one parameter). 

Every decision making act done by the Subject can be characterized by
a given boolean 3-tuple $\bold{x_1,x_2,x_3}$ of values of variables
$x_1,x_2,x_3$.
On the other hand, as pointed out above, we associate
a pure motivation $\theta_i^{\bold{x_i}}\in\Theta$ with
any $x_i,i=\overline{1,2,3}$ when $x_i=\bold{x_i}\in\{0,1\}$. 
Here 
\begin{center}
$\theta_i^{\bold{x_i}}(x_k) =
\begin{cases}
1,&  \text{ if }i<k,\\
\bold{x_i},& \text{ if } i=k,\\
0,&   \text{ if } k<i.
\end{cases}$
\end{center}

The general external order for 3-tuple $\bold{x_1,x_2,x_3}$ 
is determined by formula 
$$
\psi_{\bold{x_1,x_2,x_3}}=
\theta_1^{\bold{x_1}}\nrightarrow(\theta_2^{\bold{x_2}}\nrightarrow
\theta_3^{\bold{x_3}}).
$$

Any obtained motivation $\psi:M_*\to\{0,1\}$ determines
Subject's decision choice $\bold{x_i},i=1,2,3,$ for a given 
decision making act. 
In order to find the solution we use the following local 
extremization algorithm for 
$\psi_{\bold{x_1,x_2,x_3}}$:\\
{\bf 1) Starting at the state $x_1$ in order $(M_*,\le_*)$ 
proceed to the nearest extremum of $\theta_1^{\bold{x_1}}$.\\
2) Then continue from the found state to the nearest extremum 
of the inverted 
function $\theta_2^{\bold{x_2}}$ (due to its place in
the approximating form).\\
3) Finally, repeat starting from the found state
this time using the function $\theta_3^{\bold{x_3}}$. 
}\\
(Note, that the last step of the algorithm uses 
a double inversion of the motivation $\theta_3^{\bold{x_3}}$.) 

It is easy to check that the algorithm computes an element
of $\arg\max\psi_{\bold{x_1,x_2,x_3}}\subseteq\{x_1,x_2,x_3\}$. 
At the second step we turn a chosen $x_i$ into
$\bold{x_i}$. This number is Subject's choice in the decision making situation
defined via $\bold{x_1,x_2,x_3}$. 

It turns out that for all boolean 3-tuples $\bold{(x_1,x_2,x_3)}$
the boolean value computed with the aforementioned scheme 
coincides with the value computed with the formula $x_3\to(x_2\to x_1)$. 
Thus, Lefebvre's latter formula of behaviour can be derived from our 
model of behaviour. In our opinion, Lefebvre's subjects
are distinguished
merely by a particular internal preference: the state order $(M_*,\le_*)$.
(The binarity of the external poset is presumed in Lefebvre's problem
statement.) 

Perhaps the approach introduced in this section 
can, in principle, replace Lefebvre's axioms.
It has no need for such presumptions as "Anthropic Principle", 
"Principle of Freedom", and "Simplicity postulate". In our 
opinion, such a difference is advantageous since it appears extremely 
difficult to find a solid justification for these presumptions.

Indeed, instead of seeking a body of philosophical support
we can apply the theory
of approximating forms. Additionally, 
the restriction in this section 
by these orders does not need any special justification. 

In the continuous case the toolbox of $L$-ensembles 
not only reduces it to the Boolean case but also extends
the theory's capacity.

\section{Conclusion}

As this paper demonstrates, the classical two-valued
propositional logic can be viewed merely as a realization of 
the principle of successive approximations for 
the decision-making problem within the framework of Subject-environment 
survival. 

From this viewpoint, the classical propositional logic can 
take its beginning from the survival problem. It is also  
important that such hypothetical origin of Logics appears quite natural. 

Furthermore, this approach can serve as a background for consideration of other 
families of mappings from one poset to another with a chosen notion of 
simplicity of mapping. Any such case generates a corresponding logic.

Later in the paper we demonstrated how the effects explained 
by Lefebvre's model 
can be viewed merely as implications of 
choosing 
the binary linear external order and the three-element
linear internal preference. Such choice 
reflects on the ordinary division
of the time axis into three parts: the past, the present, and the future.

Taking into account the established connection between logic 
and approximating
forms one may say that the psychological effects described via Lefebvre's
model can be interpreted as logic of evaluation 
operators of the kind
$\psi:(\{1,2,3\},\{1<2,1<3,2<3\})\to(\{0,1\},\{0<1\})$. 
This fact can explain prevalence of the effects and partially of 
the "golden section" method.

\end{document}